\documentclass[12pt]{article}
\usepackage{latexsym}
\usepackage{amsmath,amsfonts,amssymb,theorem}
\usepackage{color}
\newcommand{\N}{\mathbb N}

\newcommand{\R}{\mathbb R}
\newcommand{\Z}{\mathbb Z}

\newcommand{\C}{\mathbb C}

\newcommand{\CA}{{\rm deg}}
\newcommand{\CB}{{\cal B}}

\newcommand{\CG}{{\cal G}}
\newcommand{\CH}{{\cal H}}

\newcommand{\CL}{{\cal L}}

\newcommand{\pR}{\R_{\geq 0}}

\newcommand{\Spec}{\operatorname{Spec}}

\newcommand{\innt}{\operatorname{int}}
\newcommand{\initial}{\operatorname{in}}
\newcommand{\init}{\initial}

\newcommand{\gHom}{\mbox{\rm Hom}}

\newcommand{\kk}{{\mathbb K}}
\newcommand{\cc}{{\mathbb C}}
\newcommand{\kx}{{\mathbf x}}
\newcommand{\Hilb}{{\rm Hilb}}
\newcommand{\km}[1]{{\mathbf #1}}
\newcommand{\HomMA}{\gHom_{\kk[\kx]}(M_X,A_X)}

\newcommand{\kitem}{\begin{itemize}\vspace{-2ex}}
\newcommand{\kenditem}{\vspace{-1ex}\end{itemize}}
\newcommand{\surj}{\rightarrow\hspace{-0.8em}\rightarrow}

\newcommand{\ko}{\overline}

\newcommand{\keps}{\varepsilon}
\newcommand{\veee}{{\scriptscriptstyle\vee}}

\newcommand{\qed}{{
\unskip\nobreak\hfil\penalty50\hskip0.1em\hbox{}\nobreak\hfil$\Box$
\parfillskip=0pt \finalhyphendemerits=0 \par\medbreak}}

\newtheorem{prop}{Proposition}
\newtheorem{lemma}[prop]{Lemma}
\newtheorem{remark}[prop]{Remark}
\newtheorem{algo}[prop]{Algorithm}
\newtheorem{ex}[prop]{Example}
\newtheorem{theorem}[prop]{Theorem}
\newtheorem{cor}[prop]{Corollary}
\newcommand{\klabel}[1]{\label{#1}\kkk{#1}}

\reversemarginpar
\newcommand{\kkk}[1]{}
\newcommand{\ktrash}[1]{{}}
\definecolor{Klaus}{rgb}{0.99,0.05,0.83}
\newcommand{\kc}[1]{#1}

\newcounter{Abschnitt}[section]

\begin{document}
\title{\bf\huge The Graph of Monomial Ideals}
\author{Klaus Altmann \and Bernd Sturmfels\footnote{Partially supported
by NSF grant DMS-0200729}}
\date{}
\maketitle
\begin{abstract}
There is a natural infinite graph whose vertices are
the monomial ideals in a polynomial ring $\kk[x_1,\dots,x_n]$.
The definition involves Gr\"obner bases or
the action of the algebraic torus $(\kk^\ast)^n$.
We present algorithms for computing
the (affine schemes representing)  edges in this graph.
  We study the induced subgraphs on
multigraded Hilbert schemes and on square-free monomial ideals.
In the latter case, the edges correspond to generalized bistellar
flips.
\end{abstract}

%
%

\section{Edge ideals}\klabel{edge}

The most important tool for computing with ideals
in a polynomial ring $\kk[\kx] = \kk[x_1,\ldots,x_n]$
over a field $\kk$ is the theory of Gr\"obner bases. It furnishes
\kc{degenerations} of arbitrary ideals in $\kk[\kx]$
to monomial ideals along one-parameter subgroups of $(\kk^*)^n$;
see \cite[\S 15.8]{Eis}.
Monomial ideals are combinatorial objects. They represent
the most special points in the ``world of ideals''.
The following adjacency relation among monomial ideals
extracts the combinatorial essence of Gr\"obner \kc{degenerations}.
\par


{\bf Definition.}
We define the infinite {\em graph of monomial ideals}
$\,\CG=\CG_{n,\kk}\,$ as follows.
The vertices of $\CG$ are the monomial ideals in $\kk[\kx]$,
and two monomial ideals $M_1, M_2 $
are connected by an edge
if there exists an ideal $I$ in $\kk[\kx]$ such that
the set of all initial monomial ideals of $I$, with respect to all
term orders,
is precisely $\{ M_1, M_2 \}$.
\par

First examples of interesting finite subgraphs
can be obtained by restricting to artinian ideals
of a fixed colength $r$.
We consider the induced subgraph on the set
\[
\CG^r \,\, := \,\, \CG^r_{n,\kk} \,\,:= \,\,
\{\, M\subseteq \kk[\kx] \mbox{ monomial ideal}\,\,\,:\,\;
\dim_\kk \kk[\kx]/M =r \, \}.
\]


\begin{prop}
\klabel{connGr}
The finite graphs $\CG^r$ are connected components of the graph $\CG$.
\end{prop}

{\sl Proof: }
Since Gr\"obner \kc{degenerations} preserve
the colength of an ideal, the graph $\CG^r$
is a union of connected components of $\CG$. Hence it
suffices  to show that $\CG^r$ is connected.
One can connect two  vertices of $\CG^r$,
i.e., two monomial ideals $M_1, M_2\subseteq \kk[\kx]$
of the same colength,
by a sequence of ``moving single boxes''
in their socles.
Hence, we may assume that the vector spaces
$M_i / (M_1 \cap M_2)$ are one-dimensional, generated by single monomials
$m_i$. But then, the ideal
\[
I \,\,\, := \,\,\,
 (M_1 \, \cap \, M_2 ) \, + \,
 \langle m_1 - m_2 \rangle
\]
provides an edge connecting $M_1$ and $M_2$ inside $\CG^r$.
\qed

The monomial ideals of colength $r$ in $\kk[x,y]$  are
in bijection with the partitions of the integer $r$. We
computed  $\CG^r_{2,\kk}$, the {\em graph of partitions},
up to $r = 13$, using the algorithm in Section \ref{algorithm}.
Here is a small example. The graph $\CG^4_{2,\kk}$
consists of five vertices and eight edges, and it equals
the cone of the vertex
$(2,2)$
over the $4$-cycle
\begin{equation}
\label{fourcycle}
\,
( 1,1,1,1 ) \, \longleftrightarrow
( 2,1,1 ) \, \longleftrightarrow
( 3,1 ) \, \longleftrightarrow
( 4 ) \, \longleftrightarrow
( 1,1,1,1).
\end{equation}
We conjecture that  $\CG^r_{2,\kk}$ is
independent of the field  $\kk$, for all $r$, but
we are still lacking a combinatorial rule
for deciding when two partitions form an edge.
\par


\begin{remark}
\klabel{infinite}
Not all connected components of the graph $\CG$ are finite.
For instance, the induced subgraph on
the principal ideals is an infinite connected component.
\end{remark}

Let us now take a closer look at the ideals which are responsible
for the edges in $\CG$. Since
monomial ideals are homogeneous with respect to the
$\Z^n$-grading of $\kk[\kx]$, one expects that edges arise
from  ideals $I$ which admit an $(n-1)$-dimensional grading.


{\bf Definition.} An ideal $I \subseteq \kk[\kx]$ is
an {\it edge providing ideal} if the set of initial monomial ideals
$\init_\prec(I)$, as $\prec$ ranges over all term orders on $\kk[\kx]$,
has cardinality two. We call $I$
an {\it edge ideal} if there exists
$c\in\Z^n$ with both positive and negative coordinates
such that $I$ is homogeneous with respect to the induced
$(\Z^n/\Z c)$-grading of $\kk[\kx]$.
\par


\begin{prop}
\klabel{edgeideal}
Every edge ideal is an edge providing ideal.
Given any edge providing ideal $I$, there
is only one non-monomial ideal $\widetilde{I}$ among its
initial ideals ${\rm in}_w(I)$, $w \in \N^n$. Moreover,
$\widetilde{I}$ is an edge ideal connecting the same vertices as $I$ does.
\end{prop}

{\sl Proof: }
The first statement holds because
\kc{generators of edge ideals have the form
$\lambda_0\kx^u+\lambda_1\kx^{u+c}+\dots +\lambda_r\kx^{u+rc}$.
Hence,
there} are only two equivalence
classes of term orders, given by $c \prec 0$ and
$c \succ 0$. For the second statement note that
the {\em Gr\"obner fan}  of $I$
is a regular polyhedral subdivision of $\R_{\geq 0}^n$
which has exactly two maximal cones. Their intersection
is an $(n-1)$-dimensional cone $C$. The unique (up to scaling)
vector $c$ perpendicular to $C$ has both positive
and negative coordinates. Fix a vector $w$ in the relative
interior of $C$.
Then $\widetilde{I}:=\init_w(I)$ is $\Z^n/\Z c$-homogeneous
and has the same two initial monomial ideals as $I$ does.
\qed

Here is an example to illustrate this for $n=2$.
The ideal $\, I = \langle x^4+ x^2 y + y^2
+ x + y + 1 \rangle \,$ is edge providing.
The unique edge ideal is
$\, \widetilde{I} \,=\,{\rm in}_{(1,2)}(I)
\, = \, \langle x^4+ x^2 y + y^2 \rangle $.



%
%

\section{Computing the graph}\klabel{algorithm}

We fix a primitive vector $c \in {\Z}^n$
with  $c_i > 0$ and $c_j < 0$ for some
 $i ,j \in \{1,2,\ldots,n\}$.
Here {\em primitive} means that the greatest common divisor of
$c_1,c_2,\ldots,c_n$ is one.


\begin{lemma}
\klabel{omegac}
For any monomial ideal $M$ in $\kk[\kx]$,
there exists an affine scheme
$\,\Omega_c(M) \,$
which
parametrizes all
$(\Z^n/\Z c)$-homogeneous ideals $I$
with $\init_{c\prec 0}(I) = M$.
\end{lemma}

{\sl Proof: }
For any minimal generator
$\kx^u$ of $M$ let $r_u$ be the largest integer
such that $\,u + r_u\, c\, $ is non-negative. Introduce unknown
coefficients $\, \lambda_{u,1},\ldots, \lambda_{u,r_u}\,$
and form
\begin{equation}
\label{algorithm-Omega-equ}
\underline{ \kx^u }\, + \,
 \lambda_{u,1} \,\kx^{u + c} \, + \,
 \lambda_{u,2} \,\kx^{u + 2 c}  \, + \,\cdots \, + \,
 \lambda_{u,r_u} \,\kx^{u + r_u c} .
\end{equation}
The ideal $I$ generated by the polynomials (\ref{algorithm-Omega-equ})
satisfies  $\init_{c\prec 0}(I) = M$
if and only if they form a Gr\"obner basis
with the underlined leading terms.
By Buchberger's criterion, this means that all S-pairs reduce to zero,
giving an explicit system of polynomial equations in terms of the
$ \lambda_{u,i}$.
On the other hand, we would like the coordinates $ \lambda_{u,i}$
to be uniquely determined from $I$.
This is the case if we require
that (\ref{algorithm-Omega-equ}) describes
 a {\em reduced} Gr\"obner basis,
imposing $\lambda_{u,i}=0$ whenever $\kx^{u+ic}\in M$.
\qed

We call  $\,\Omega_c(M) \,$ the {\em Schubert scheme}
of $M$ in direction $c$. In the case when $M$
is generated by a subset of the variables
then $\Omega_c(M)$ is a
Schubert cell in the Grassmannian.
If $M_1,M_2$ are two monomial ideals, then
the scheme-theoretic intersection
\vspace{-1ex} \[
\Omega_c(M_1,M_2) \quad := \quad \Omega_{c}(M_1)\,\cap\, \Omega_{-c}(M_2)
\]
parametrizes all {\em $c$-edge ideals}  between
$M_1$ and $M_2$.


\begin{algo}
\klabel{algoc}
{\rm (Input: $\, c,M_1,M_2$. Output: $\,\Omega_c(M_1,M_2)$)}
\vspace{-1ex}
\end{algo}

\noindent {\sl Step 1: }
Construct the affine scheme
$\Omega_{c}(M_1)\,$ using the procedure in the proof above.
Using $S$-pair reduction,
one obtains a set of polynomials in variables $\lambda_{u,i}$, and
the universal $c$-edge ideal over
the base $\Omega_{c}(M_1)\,$ is described
by the polynomials (\ref{algorithm-Omega-equ}).
\par

\noindent {\sl Step 2: }
Construct the affine scheme
$\,\Omega_{-c}(M_2)\,$ as in Step 1. This gives
a set of polynomials in some other variables
$\widetilde \lambda_{\tilde u,i}$ representing the
universal $c$-edge ideal.
\par

\noindent {\sl Step 3: }
Form additional joint equations
in both sets of variables $\lambda_{u,i}$ and
$\widetilde \lambda_{\tilde u,i}$ which
express the requirement that the universal ideal over
$\Omega_{c}(M_1)\,$ coincides with the universal ideal
over $\Omega_{-c}(M_2)$. This is done by reducing the polynomials
(\ref{algorithm-Omega-equ}) of Step~1 modulo those of Step 2
and reading off the coefficients with respect to $\kx$.
\par

Let us demonstrate how Algorithm \ref{algoc}
works for a small example.
\par


\begin{ex}
\klabel{exalgoc}
\rm
Let
$M_1 \, = \, \langle x^6, x^2 y , y^2 \rangle$,
$\,M_2 \, = \, \langle x^2, x y^2 , y^6 \rangle\,$
and  $c=(1,-1)$.
In Step 1 we introduce three indeterminates $a_1,a_2,a_3$.
The ideals in $\Omega_{c}(M_1)$ are of the form
\begin{equation}
\label{algorithm-example-idealv}
\langle\, \underline{ x^6 } \, , \,\,
\underline{ x^2 y } \, + \, a_1 x^3 \, , \,\,
\underline{ y^2 } \, + \, a_2 xy \, +\, a_3 x^2 \rangle .
\end{equation}
These polynomials are a Gr\"obner basis
with underlined leading terms if and only if
\begin{equation}
\label{algorithm-example-basev}
a_1^2 - a_1 a_2 + a_3 \quad = \quad 0\,.
\end{equation}
In Step 2 we similarly compute the affine scheme
$\Omega_{-c}(M_2)$ to be the hypersurface
\begin{equation}
\label{algorithm-example-base-v}
b_3^2 - b_1 b_3 + b_2 \quad = \quad 0 \,,
\end{equation}
carrying the universal ideal
\begin{equation}
\label{algorithm-example-ideal-v}
\langle \,  \underline{x^2}
\, +\, b_1 x y
\, + \, b_2 y^2 \,, \,\,
\underline{x  y^2} \, +\,b_3 y^3 \,, \,\,
\underline{ y^6} \, \rangle.
\end{equation}
Finally, in Step 3 we enforce the condition that the
ideals in (\ref{algorithm-example-idealv}) and
(\ref{algorithm-example-ideal-v}) are equal, given
that (\ref{algorithm-example-basev}) and
(\ref{algorithm-example-base-v}) hold. This is done by
reducing the generators of (\ref{algorithm-example-idealv}) modulo
the Gr\"obner basis (\ref{algorithm-example-ideal-v}) and collecting
coefficients in the normal forms. We obtain
\vspace{-1ex}
\begin{equation}
\label{algorithm-example-ideal-v2}
\bigl\{ \,
a_1 - a_3 b_1 + a_3 b_3 \, , \,\,
a_2 - a_3 b_1 \, , \,\,
a_3 b_2 - 1 \, ,  \,\,
b_2-b_1 b_3 + b_3^2 \, \bigr\}.
\end{equation}
\end{ex}


\begin{ex}
\klabel{kepsilon}
\rm
The Schubert schemes $\Omega_c(M_i)$ in the previous example
are  reduced and irreducible. However, this is not true in general.
For instance, for
 $\,M=\langle
x^6,\, y^5,\, z^9,\, y^3z^5,\, x^4y^3z^2,\, x^3y^2z^4,\, x^2y^4z^3
\rangle \,$
we obtain $\,\Omega_{(-3, 0, 1)}(M) \, \simeq \,
\Spec \kk[\keps]/(\keps^2)$.
\end{ex}

Our next result will imply that the lower index ``$c$''
can be dropped from $\Omega_c(M_1,M_2)$.


\begin{theorem}
\klabel{only1c}
Given any two monomial ideals $M_1,M_2$ in  $\kk[\kx]$,
there is at most one direction $\,c \in \Z^n \,$
such that the scheme $\,\Omega_c(M_1,M_2)\,$ is non-empty.
Moreover, if $\,\Omega_c(M_1,M_2)\neq\emptyset\,$,
then
$M_1,M_2$ have equal Hilbert functions with respect to an
induced $(\Z^n/\Z c')$-grading if and only if $c'=\pm c$.
\end{theorem}

The proof of Theorem \ref{only1c} will be given in the next section.
If $M_1$ and $M_2$  are connected by an edge
in our graph $\CG$, then $c$ is uniquely determined, and we simply write
$$ \Omega(M_1,M_2) \quad := \quad \Omega_c(M_1,M_2) $$
for the scheme which parameterizes all edge
ideals between $M_1 $ and $M_2$.
If $M_1$ and $M_2$  are not connected by an edge
in $\CG$ then $ \Omega(M_1,M_2) $ denotes the empty set.
Hence the following algorithm can be used to
determine the adjacency relation in $\CG$.


\begin{algo}
\klabel{algofact}
\rm
(Input: $\, M_1,M_2$. Output: $\,\Omega(M_1,M_2)$)
\vspace{-1ex}
\end{algo}

\noindent {\sl Step 1: }
Compute the $\N^n$-graded Hilbert series
$H(M_i; \kx)$ of the two given monomial ideals
as rational functions, i.e., find the
numerator polynomials $K_1$ and $K_2$ of
$$ H(M_i; \kx) \quad = \quad \frac{K_i(x_1,\ldots,x_n)}{
 (1-x_1) (1-x_2) \cdots (1-x_n)}
$$

\noindent {\sl Step 2: }
Factor the polynomial $\,K_1(\kx) -  K_2(\kx)\,$
into irreducible factors.
Output $\Omega_c(M_1,M_2)=\emptyset$,
unless there
\kc{is, up to sign, a unique}  primitive vector
$c \in \Z^n$
which has positive and negative coordinates
such that
the binomial $\kx^{c_+} -  \kx^{c_-}$
appears as a factor.

\noindent {\sl Step 3: }
Run Algorithm \ref{algoc} for
the vector $\pm c$ found in Step 2,
and output the affine scheme
$ \,\Omega(M_1,M_2) \, = \, \Omega_c(M_1,M_2)$.
(It is still possible that this scheme empty.)

\smallskip

The correctness of Algorithm \ref{algofact} follows
directly from Theorem \ref{only1c}.
An improvement to Step 1 in this algorithm
in the context of an ambient gradient
will be discussed in the next section.

As an example consider the two ideals in Example \ref{exalgoc}.
In step 1 we compute
$$ K_1(x,y) =  1 -  x^6  -  x^2 y -  y^2 +
 x^6  y  + x^2  y^2 $$
$$ K_2(x,y) =  1 -  x^2 - x y^2 - y^6  +   x^2 y^2 + x y^6  .$$
The difference  $\, K_1(x,y) - K_2(x,y) \,$
of these numerator polynomials factors as
$$
(x-y)(y-1)(x-1)(x^4+x^3y+x^2y^2+xy^3+y^4+x^3+x^2y+xy^2+y^3+x^2+xy+y^2+x+y).
$$
The only binomial factor with both terms non-constant is
$\,x-y $, and we conclude that
$\,\Omega(M_1,M_2) \,$ equals $\,\Omega_{(1,-1)}(M_1,M_2) $,
the affine scheme described by (\ref{algorithm-example-ideal-v2}).

%
%

\section{Multigraded Hilbert schemes}\klabel{Hilbert}

We consider an arbitrary grading of the polynomial ring
$\kk[\kx]$. It is given by an epimorphism of abelian groups
$ \, {\rm deg} \, : \, \Z^n \surj A$.
For any function $\,h :A \rightarrow \N$, the
{\it multigraded Hilbert scheme} $\,\Hilb_h \,$
parametrizes all homogeneous ideals $I$ such that
$\kk[\kx]/I$ has Hilbert function $h$.
This scheme was introduced in \cite{HaSt}.
Multigraded Hilbert schemes provide a natural setting
for studying finite subgraphs of $\CG={\CG}_{n,\kk}$.
\par


{\bf Definition.}
A multigraded Hilbert scheme $\Hilb_h$
has the {\em induced subgraph property} if any two monomial ideals
$M_1, M_2\in\Hilb_h$ which are connected in ${\cal G}_{n,k}$
can also be connected via an edge ideal $I$  which lies in
the same Hilbert scheme $\Hilb_h$.
\par

The induced subgraph property holds for the Hilbert scheme
of points, where $A = \{0\}$ is the zero group,
by our discussion in Section \ref{edge}.
 However, it fails in general.
\par


\begin{ex}
\klabel{super}
\rm
Consider the ``super-grading'' of  $\kk[x,y]$ given by
${\rm deg} :\Z^2\to \Z/2\Z,$ $ \,(r,s)\mapsto r+s$,
and define $ h : \Z/2\Z \rightarrow \N$ by
$h(0) = h(1) = 2$.  The two ideals
$\,M_1 = \langle x^4, y \rangle \,$ and
$\,M_2 = \langle x^2, y^2 \rangle \,$
are points in  $\Hilb_h$.
They are connected in $\CG$ as was seen in
(\ref{fourcycle}). Algorithm \ref{algofact} finds
that the edge ideals are
$\, \langle x^2+ \alpha  y \, ,\, y^2 \rangle \,$
for any $\alpha \in \kk^\ast$.
None of the edge ideals is
homogeneous in the given grading.
We conclude that the Hilbert scheme $\,\Hilb_h \,$ does not
have the induced subgraph property.
\end{ex}


{\bf Definition. }
A grading of $\kk[\kx]$ is called {\em positive} if only the constants
have degree $0$. This implies that  the grading group $A$ is torsion-free,
i.e., $A\cong \Z^q$ for some $q$.
\par

A torsion-free grading $\CA:\Z^n\to \Z^q$ is positive if and only if
$\,\N^n\cap \ker(\CA)=0 \,$ if and only if the
{\em fibers} $\,\N^n\cap\CA^{-1}(a)$ are finite
if and only if the polyhedra
$\,\pR^n\cap\CA_{\R}^{-1}(a)\,$ are compact.
Under these circumstances, our graphs behave nicely:
\par


\begin{theorem}
\klabel{inducedsub}
Let $\,\CA \!:\Z^n\to \Z^q$ be a positive grading and
$h :\Z^q\to\N$ any function.
Then the multigraded Hilbert scheme
$\,\Hilb_h$ has the induced subgraph property.
\end{theorem}

We will derive this theorem from the following lemma.


\begin{lemma}
\klabel{degc0}
Let $\CA \! :\Z^n\to \Z^q$ be a positive grading and
$M_1,M_2 \subset  \kk[\kx]$ monomial ideals
with the same Hilbert function.
Then  $\, \Omega_c(M_1,M_2) \not= \emptyset \,$ implies
$\CA(c)=0$.
\end{lemma}

{\sl Proof of Theorem \ref{inducedsub}: }
Let $M_1$ and $M_2$ be monomial ideals in $\Hilb_h$
and $I $ an edge ideal in $\Omega(M_1,M_2)$.
Lemma \ref{degc0} implies that $I$ is homogeneous with respect
to the given positive grading $\, {\rm deg}$.
Since $M_1$ and $M_2$ are initial ideals of $I$,
all three ideals have the same Hilbert function,
and hence $I$ is a point in  $\Hilb_h$ as desired.
\qed

{\sl Proof of Lemma \ref{degc0}: }
Let $I \in \Omega_c(M_1,M_2)$,
$M_1=\initial_{c\prec 0}(I)$ and $M_2=\initial_{c\succ 0}(I)$.
The edge ideal  $I$ is generated by $\Z^n/\Z c$-homogeneous polynomials
of the form
\begin{equation}
\label{chomopoly}
  { \kx}^u \, + \,
 \lambda_1 { \kx}^{u + c} \, + \,
 \lambda_2 { \kx}^{u + 2 c}  \, + \,\cdots \, + \,
 \lambda_r { \kx}^{u + r c} \quad \qquad (\lambda_r\neq 0).
\end{equation}
We shall abuse the symbols $M_1,M_2,I$
to also denote the set of exponents of the
monomials in that ideal. For instance,
from (\ref{chomopoly}) we infer $u\in M_1$ and
$u+rc\in M_2$. We also have the following obvious inclusions
among finite sets of monomials:
\begin{equation}
\label{obviousinclusion}
 I \cap \CA^{-1}(a) \,\subseteq \, M_i\cap \CA^{-1}(a)
\qquad \hbox{for $i=1,2$
and $a \in \Z^q$}.
\end{equation}
Our strategy is this: we first prove Lemma \ref{degc0} for one-dimensional
gradings.

{\sl Step 1: $\;q=1$.}
Assume that $\,c \not\in \ker(\CA)$.
We claim that
\begin{equation}
\label{notinclusion}
M_1\cap \CA^{-1}(a) \;\subseteq\; M_2\cap \CA^{-1}(a)
\;\subseteq\; I\cap \CA^{-1}(a)
\qquad \hbox{for all $a \in \Z^q$}.
\end{equation}
This implies
$M_1\cap \CA^{-1}(a) =M_2\cap \CA^{-1}(a)$, hence
$ M_1 = M_2 $, a
contradiction  which will establish Lemma 13 for $q=1$.

We may assume  $\CA(\N^n)\subseteq \N$ and $\CA(c)<0$.
We shall prove (\ref{notinclusion}) for positive
integers $a$ by induction. The case $a \leq 0$ is void.
Suppose the two inclusions hold for all $a<a_0$.
Consider any element $\,u\in M_1\cap \CA^{-1}(a_0)\,$
and a corresponding polynomial
$\,f= \kx^u + \lambda_1 \kx^{u+c}  +  \cdots  +
\lambda_r \kx^{u+rc} \in I$ with
$\lambda_r\neq 0$ and minimal $r\geq 0$.
If $r=0$, then $\kx^u\in I$, hence $u\in M_2$.
If $r > 0$, then $u+rc\in M_2$ with $\CA(u+rc)=a_0+r\cdot\CA(c)<a_0$.
This  implies $u+rc\in I$ by the induction hypothesis. But then $f$
can be shortened, and we obtain a contradiction.
The Claim  (\ref{notinclusion}) follows.

{\sl Step 2: $\;q\geq 2$.}
Consider the polyhedral cone
$\sigma:=\CA_\R(\pR^n)$ in $\R^q$.
Since $\N^n\cap \ker (\CA)=0$, the cone
$\sigma$ is pointed which means  that
the dual cone $\sigma^\vee$ is full-dimensional.
For a linear map  $\ell:\Z^q\to\Z$ the following
statements are equivalent:
\vspace{-0.8ex}
\[
\renewcommand{\arraystretch}{1.2}
\begin{array}{rcl}
\N^n\cap \ker (\ell\circ\CA)=0
&\Longleftrightarrow&
\N^n\cap\CA^{-1}\big(\ell^{-1}(a)\big) \mbox{ are finite for all $a\in\Z$}\\
&\Longleftrightarrow&
\CA(\N^n) \cap \ell^{-1}(a) = \sigma \cap \ell^{-1}(a)\,
	\mbox{ are finite}\\
&\Longleftrightarrow&
\sigma\cap (\ker \ell) = 0\\
&\Longleftrightarrow&
\ell\in (\innt \sigma^\veee) \cup (-  \innt\sigma^\veee)\,
\vspace{-0.8ex}
\end{array}
\]
Fix a basis $\CB$ of $(\R^q)^\ast$
consisting of linear forms $\ell$ which satisfy this condition.
For each $\ell \in \CB$,
we apply Step 1 to the one-dimensional grading
$\,(\ell\circ\CA):\Z^n\to\Z^q\to\Z$, and we conclude that
$c$ lies in  $\ker (\ell\circ\CA)$. Therefore,
$\,c \in \bigcap_{\ell \in \CB} \ker (\ell \circ\CA) =\ker( \CA) $,
since $\R \CB =  (\R^q)^\ast$.
This finishes the proof of Lemma \ref{degc0} and of
Theorem~\ref{inducedsub}.
\qed

Suppose that $M_1$ and $M_2$ are monomial ideals
on a multigraded Hilbert scheme $\Hilb_h$.
For $a\in A$ we denote by $\,P_a(M_i) \in \N^n$
the sum of all vectors
$u \in \N^n$ such that $\kx^u\notin M_i$
and $\,\CA(u)=a$. Here the number of summands
 is $h(a)$, the value of
the Hilbert function at $a$.


\begin{lemma}
\klabel{posmultiple}
Let $M_1,M_2 \in \Hilb_h$,
$\Omega_c(M_1,M_2) \not=\emptyset$, and ${\rm deg}(c) = 0$.
If $M_1,M_2$ differ in a degree  $a\in A$, then
$\,P_a(M_1) - P_a(M_2)$
are positive integer multiples of $c$.
\end{lemma}

{\sl Proof: }
Let $I\in\Omega_c(M_1,M_2)$.
We may assume that  $\CA$ equals
the $c$-grading $\,\Z^n\to \Z^n/\Z c$.
For a degree $a\in \Z^n/\Z c$ we denote by $I_a$ and $(M_i)_a$ the
homogeneous parts of the corresponding ideals.
Let ${\cal L}$ be a finite set of polynomials such that
$(M_1)_a$ and $(M_2)_a$ are contained in $\init_{c\prec 0}(\CL)$
and $\init_{c\succ 0}(\CL)$, respectively.
For an element $\lambda_0\,{\bf x}^u+\dots+\lambda_r\,
{\bf x}^{u+rc}\in\CL$
with $\lambda_0,\lambda_r\neq 0$ we call $r$ its length.
The {\em total length} of $\CL$ is the sum of the lengths
of all polynomials in $\CL$.
Now, whenever there are two elements $f,g\in\CL$ having the same highest
or the same lowest monomial, then we can reduce the total length
of ${\cal L}$ without
loosing $(M_1)_a\subseteq\init_{c\prec 0}\CL$
and $(M_2)_a\subseteq\init_{c\succ 0}\CL$.
Just replace $\{f,g\}$ by the shorter polynomial among them and $ f-g$.
Iterating this several times, we arrive at a set
$\CL$ none of whose polynomials have common ends.
The set $\CL$ provides a bijection
$(M_1)_a \stackrel{\sim}{\rightarrow} (M_2)_a$ via
$\,\init_{c\prec 0}(f)\mapsto \init_{c\succ 0}(f)$.
\qed

We are now prepared to tie up some loose ends from the last section.
Let us first reexamine the  process
of finding the correct direction $c$ in
Algorithm \ref{algofact}. Factoring the
numerator difference of the Hilbert series
can be replaced by the following procedure.
\par


\begin{algo}
\klabel{algohilb}
\rm
(Input: $\, M_1,M_2 \in \Hilb_h$ with respect to a
positive grading or $A=0$. Output: $\,\Omega(M_1,M_2)$)
\vspace{-1ex}
\end{algo}

\noindent {\sl Step 1: }
Pick a degree $a \in A$ in which the monomial ideals $M_1$ and $M_2$
are different. Compute the vectors $P_a(M_1)$ and $P_a(M_2)$.

\noindent {\sl Step 2: }
If $P_a(M_1) = P_a(M_2)$ then stop and output the empty set.
Otherwise let $c$ be the primitive vector
in direction  $P_a(M_1)- P_a(M_2)$.

\noindent {\sl Step 3: }
 Using Algorithm \ref{algoc}, compute and output $\Omega_c(M_1,M_2)$.

\smallskip

Finally, it is time to present the

{\sl Proof of Theorem \ref{only1c}: }
Suppose that $\,\Omega_c(M_1,M_2)\,$
and $\,\Omega_{c'}(M_1,M_2)\,$ are both non-empty,
where $c$ and $c'$ are primitive vectors in $\Z^n$
which have positive and negative coordinates.
The group $\,A := \Z^n/\Z {c'} \simeq \Z^{n-1}\,$
is torsion-free and the canonical map
$ \, {\rm deg} \, : \, \Z^n \surj A\,$
is a positive grading. Applying Lemma \ref{degc0}
to this grading, we find that that $c=\pm c'$.
Finally, Lemma \ref{posmultiple} excludes $c'=-c$.
\qed

One important question regarding Hilbert schemes is
under which circumstances $\Hilb_h$ is connected. While classical
Hilbert schemes are known to be
connected \cite{Ree}, Santos \cite{San} recently constructed a
disconnected multigraded Hilbert scheme.
The graph  introduced in
this paper provides a tool for studying connectivity questions.


{\bf Definition.}
For a subscheme $\CH\subseteq \Hilb_h$,
we denote by $\CG(\CH)\subseteq\CG$ the subgraph with vertices and edges
built from monomial and edge ideals in $\CH$.
In particular, the induced subgraph property means that
$\CG(\Hilb_h)$ is an induced subgraph of $\CG$.
\par


\begin{lemma}
\klabel{connGH}
Let $\CA \!:\Z^n\to\Z^q$ be a positive grading of $\kk[\kx]$
where $\kk = \R$ or $\kk = \cc$.
If $\CH $ is an irreducible component of $\Hilb_h$
then the graph $\CG(\CH)$ is connected.
\end{lemma}

{\sl Proof: }
The positive grading implies that $\Hilb_h$ is a projective scheme
\cite[Corollary 1.2]{HaSt}. Hence $\CH$ is irreducible and projective.
The algebraic torus $(\kk^*)^n$ acts on $\CH$ with finitely many
fixed points (the monomial ideals).
Consider any two monomial ideals $M_1,M_2$ which lie in $\CH$.
Then $\{M_1,M_2\}$ is  an edge in $\CG(\CH)$ if and only if
$M_1$ and $M_2$ are in
the closure of a one-dimensional torus orbit on $\CH$.
The irreducible variety $\CH$ contains a
 connected projective curve $C$, not necessarily
irreducible, which lies in $\CH$ and contains both points
$M_1 $ and $M_2$. We can \kc{degenerate} the curve $C$ by
a generic one-parameter subgroup of $(\kk^*)^n$ to
a curve $C'$ which is $(\kk^*)^n$-invariant. This can be done,
for instance, by a Gr\"obner basis computation in the
homogeneous coordinates of the projective variety $\CH$.
The \kc{degenerate} curve $C'$ still contains $M_1$ and $M_2$,
it is connected \kc{(since, by
Stein Factorization, flat degenerations
of connected projective schemes are connected;
see e.g.~Exercise III/11.4 in \cite{Ha}),}
and it is set-theoretically a union
of closures of one-dimensional torus orbit on $\CH$.
Hence $M_1$ can be connected to $M_2$ by a sequence
of edges in $\CG(\CH)$.
\qed

\begin{cor}
\klabel{connHilb}
For positive gradings with $\kk = \R$ or $\kk = \cc$,
the multigraded Hilbert scheme $\Hilb_h$ is connected if
and only if the graph ${\CG}(\Hilb_h)$ is connected.
\end{cor}

{\sl Proof: }
The if direction always holds even if the grading is not positive
and  $\Hilb_h$ is not compact. Indeed,  if $I_1$ and $I_2$
are arbitrary ideals in   $\Hilb_h$ then we can connect
them to their initial ideals ${\rm in}_\prec(I_1)$
and ${\rm in}_\prec(I_2)$  under some term order $\prec$.
Connecting these two monomials ideals
along the graph ${\CG}(\Hilb_h)$ establishes
a path in $\Hilb_h$ which connects $I_1$ and $I_2$.
For the only-if direction we use Lemma \ref{connGH}.
Suppose $\Hilb_h$ is connected.  Then the graph
of irreducible components is connected, where two
components are connected by an edge in this graph
if and only if they intersect.
On the other hand, with $\Hilb_h$, all its irreducible components
are torus invariant. Hence,
by Gr\"obner \kc{degenerations},
every non-empty intersection of irreducible components of $\Hilb_h$
contains at least one monomial ideal. Using Lemma \ref{connGH},
we can then connect any two
monomial ideals $M_1,M_2 \in \Hilb_h$ by
a sequence of edges in $\CG(\Hilb_h)$.
\qed

We do not know at present whether Lemma \ref{connGH} and
Corollary \ref{connHilb}
remain valid if the grading is not positive.
Corollary \ref{connHilb} had been proved previously by
Maclagan and Thomas for the special case
of {\em toric Hilbert schemes} \cite{MaTh}.
Here ``toric'' means that $h$ is the characteristic
function of ${\rm deg}(\N^n)$.
The disconnected example in \cite{San} is a toric
Hilbert scheme. It was constructed using methods
from polyhedral geometry.

%
%

\section{Simplicial complexes}\klabel{simplicial}

Every class of monomial ideals determines
an induced subgraph of $\CG$.
In this section we study the induced finite subgraph
on square-free monomial ideals in $\kk[\kx]$.
These ideals correspond to simplicial complexes on   $[n]:=\{1,2,\dots,n\}$.
We write $\Delta_{n-1 }$ for the full simplex on $[n]$.
Faces of $\Delta_{n-1}$ are subsets of $[n]$,
and they are identified with their
incidence vectors in $\{0,1\}^n$.
Fix an arbitrary simplicial complex $X \subset \Delta_{n-1}$.
Its {\em Stanley-Reisner ideal}
and  its {\em Stanley-Reisner ring} are
\vspace{-1ex}
\[
M_X \,\,:= \,\,\langle \,\kx^u \,\,:\,\, u\in \Delta_{n-1}\backslash X\,
\rangle \;\subseteq\; \kk[\kx]
\hspace{1em}\mbox{and}\hspace{1em}
A_X \,\, := \,\, \kk[\kx]/M_X\,.
\vspace{-0.5ex}
\]
The  $A_X$-module $\,\HomMA\,$
describes the {\em infinitesimal deformations} of $A_X$.
It is $\Z^n$-graded. Elements $\lambda$
of degree $c$ in $\,\HomMA\,$
look like $\,\kx^u\mapsto \lambda(u)\, \kx^{u+c} $,
where $\lambda$ ranges over a subspace
of the vector space of maps
$\Delta_{n-1}\backslash X \to \kk$,
\kc{cf.~\cite{AlCh}}.
The equations defining this subspace include
$\, \lambda(u) = 0 \,$ whenever
$\,u+c \not\in \N^n $.
For any $c \in \Z^n$ and any $\lambda \in
\HomMA_c $, we define an ideal as follows:
$$ I_\lambda \quad := \quad
\langle \,\,\kx^u + \lambda(u)\, \kx^{u+c} \, : \,
u \in \Delta_{n-1}\backslash X \,\rangle. $$
If $\kx^{u+c}\in M_X$, then the value $\lambda(u)$ does not matter
neither for $\lambda\in\HomMA$, nor for $I_\lambda$.
We will set $\lambda(u):=0$ in this case,
cf.~the end of the proof of Lemma \ref{omegac}.
\par


\begin{theorem}
\klabel{T1Omega}
Let $c\in \Z^n$ be a vector with both positive and negative coordinates.
\vspace{0.5ex}\\
(a) The map $\HomMA_c \rightarrow \Omega_c(M_X),
 \lambda \mapsto I_\lambda $ is an isomorphism of schemes
over $\kk$.
In particular, the Schubert scheme $\,\Omega_c(M_X)\,$ is an affine space.
\vspace{0.5ex}\\
(b) The monomial $c$-neighbors of $M_X$ in $\CG$
come from $\,\HomMA_c$ via
\vspace{-1ex}
\[
M_X^\prime(\lambda)\, :=\,
\langle \kx^u : \, u\in
\Delta_{n-1}\backslash X,\; \lambda(u)=0  \rangle
\; +\;
\langle  \kx^{v+c}:\, v\in
\Delta_{n-1}\backslash X,\; \lambda(v)\neq 0 \rangle.
\vspace{-1ex}
\]
\end{theorem}

{\sl Proof: }
(a)  Each pair $(\kx^u,\kx^v)$ of minimal $M_X$-generators
provides a  condition on both sides, in addition
to the previously mentioned vanishing of
certain $\lambda$-coordinates. The condition
is gotten via the linearity of $\lambda\in\gHom$, on the one hand, and
via the S-polynomials, on the other.
In both cases, one obtains that
$\lambda(u)=\lambda(v)$ whenever $\kx^{(u\cup v)+c}\notin M_X$.
In particular, these equations are linear.
\vspace{0.5ex}\\
(b) We must show that the generators
$\kx^u + \lambda(u) \kx^{u+c}$
with $u\in \Delta_{n-1} \backslash X$
form a Gr\"obner basis
of $I_\lambda$ also for the term order
$\,c\succ 0$.
Let $\kx^v=\init_{c\succ 0} (f)$ be the
 initial term of some element $f \in I_\lambda$.
We must show that $\kx^v$ is a multiple of
the $(c\succ 0)$-leading term of some $\kx^u + \lambda(u) \kx^{u+c}$.
After reducing $f$ to normal form with respect
to the generators, only two cases remain. Either
$f$ is a binomial or a monomial.
\vspace{0.5ex}\\
{\sl Case 1: $f$ equals $\,\kx^{v-c-u}(\kx^u + \lambda(u) \kx^{u+c})$
with $\lambda(u)\neq 0$.}
Then $\kx^v$ is divisible by
$\kx^{u+c}=\init_{c\succ 0}(\kx^u + \lambda(u) \kx^{u+c})$, and we are done.
\vspace{0.5ex}\\
{\sl Case 2: $f$ equals $\kx^v$, i.e., $\kx^v\in I_\lambda$.}
For $w\in\N^n$ let $\ko{w} = \{i : w_i \not= 0 \}$ denote
its support. Then $\ko{w_1+w_2}=\ko{w_1}\cup\ko{w_2}$, and
$\ko{w}=w$ for elements $w\in \Delta_{n-1}$.
 The ideal $M_X$ being square-free, we have
$\kx^w\in M_X$ if and only if $ \kx^{\ko{w}}\in M_X$.
In particular, since $\kx^v\in {\rm in}_{c \prec 0}(  I_\lambda)
= M_X$, we have $ \kx^{\ko{v}}\in M_X$ and
$\ko{v} \in \Delta_{n-1} \backslash X$.
It suffices to show that $\lambda(\ko{v})=0$.
Suppose $\lambda(\ko{v})\not=0$.
Then  $v+c\geq \ko{v}+c\geq 0$.
Now,
$\lambda(\ko{v})\,\kx^{v+c}=\kx^{v-\ko{v}}(\kx^{\ko{v}} +
\lambda(\ko{v}) \kx^{\ko{v}+c}) - \kx^v$
implies that $\kx^{v+c}\in I_\lambda$,
hence $\ko{v+c}\in\Delta_{n-1}\backslash X$.
Setting $w_1:=\ko{v}$ and $w_2:=\ko{v+c}$,
we find $\ko{(w_1\cup w_2)+c}=\ko{\ko{v}+c}\,\in\, X$
(since $\lambda(\ko{v})\neq 0$), i.e.,
$\kx^{(w_1\cup w_2)+c}\notin M_X$.
The equations mentioned in (a) imply
$\lambda(\ko{v})=\lambda(\ko{v+c})$.
We can now replace by $v$ by $v+c$
and run the same argument again. After
iterating this step finitely many times,
the hypothesis $v+c \geq 0$ will no longer hold,
so that $\lambda(\ko{v+c}) = 0$ and hence
$\lambda(\ko{v}) = 0$. This completes the proof.
\qed

To make the previous theorem more useful, we shall apply the
description of the vector spaces
$\,\HomMA_c \,$ given by
Altmann and Christophersen in \cite{AlCh}:

{\bf Notation.}
For a subset $N\subseteq X$,
we denote by
$\langle N \rangle$ the union
of all {\em open} simplices $|f|$,
$f\in N$,
in the geometric realization $|X|$.
For $c\in \Z^n$ with non-trivial positive and negative parts
$c^+$ and $c^-$, we denote by $a,b\subseteq[n]$ their
respective supports, and
\[
\renewcommand{\arraystretch}{1.2}
\begin{array}{rcl}
N_c&:=& \big\{ \,f\in X \,\, \, : \,\,\,
a\subseteq f\,,\; f\cap b =\emptyset\,,\; f\cup b \notin X\big\}\,,\\
\tilde{N}_c&:=& \big\{ f\in N_c \,\,\,:\,\,\,
 f\cup b^\prime \notin X
\,\,\,\hbox{for some proper subset $b'$ of $b$} \,\big\}.
\end{array}
\vspace{-2ex}
\]
\par

The following results are proved in  \cite{AlCh}.
If $c_i \leq -2$ for some $i$ then
$\HomMA_c$ \kc{vanishs}. If not, i.e.,
if $c^-=b$, let  $\,N^1,\dots, N^m$ be the subsets of $N_c$
which correspond to those connected components
of $\langle N_c \rangle$ that do not touch $\tilde{N}_c$.
There is an  isomorphism
$\,\kk^m\stackrel{\sim}{\rightarrow}\HomMA_c$.
It sends  $(\lambda_1,\dots,\lambda_m)$ to the map
$\lambda: \Delta_{n-1}\backslash X \to \kk$ defined as
$\,\lambda(u):=\lambda_i\,$ if $\,(u\cup a)\backslash b \in N^i$
and $\,\lambda(u):=0$ otherwise.
Theorem \ref{T1Omega} (a) implies that
$\Omega_c(M_X)$ is trivial
unless $\,a = {\rm supp}(c^+)\,$ is a face of $X$
\kc{($a\notin X$ $\Rightarrow$ $N_c=\emptyset$)}.

Suppose $a = {\rm supp}(c^+) \in X$ and $c^- = b$ and
fix $N^1,\ldots,N^m$ as above. Then each non-empty subset
 $\{i_1,\dots,i_\ell\}\subseteq\{1,\dots,m\}$
determines a monomial ideal as follows:
\begin{eqnarray*}
M_X^\prime \;=\; M_X^\prime(i_1,\dots,i_\ell)\;=\;
& \langle \,\, \kx^u \,\,\,:\,\,\, u\in
\Delta_{n-1}\backslash X\,,\;
(u\cup a)\backslash b \,\notin N^{i_1}\cup\dots\cup N^{i_\ell}
\,\, \rangle \\
+ \,\, & \! \langle \,
\kx^{v+c}\,\,:\,\; v\in
\Delta_{n-1}\backslash X \,,\;
(v\cup a)\backslash b \,\in N^{i_1}\cup\dots\cup N^{i_\ell} \,\, \rangle .
\end{eqnarray*}
These $2^m-1$ ideals are generally not distinct..
However, if $(a\cup b)\notin X$,
then $\langle N_c \rangle$ is connected, hence
$m=1$ for $\tilde{N}_c=\emptyset$ and
$m=0$ for $\tilde{N}_c\neq\emptyset$.
Theorem \ref{T1Omega} and the results quoted from \cite{AlCh} imply


\begin{cor}
\klabel{neighbors}
The ideals $M_X^\prime$
are all the neighbors of $M_X$ in $\CG$ in direction $c$.
\end{cor}

We next identify  the square-free monomial ideals among
the neighbors $M_X^\prime$ of $M_X$. From the generators
$\,\kx^{v+c}\,$ we see that $M^\prime_X$ is not square-free
unless $c_i \leq 1 $ for all $i$. Hence from now on we assume
that $c^+=a$ and $c^-=b$ are non-empty disjoint subsets of $[n]$,
\kc{and, w.l.o.g., $\{N^{1},\dots,N^{\ell}\}$}
is a non-empty subset of the
connected components of
$\langle N_c \rangle$ that do not touch $\tilde{N}_c$.
With these data we associate the distinguished subcomplex
\[ F \quad := \quad \big\{ f\backslash a
\,\,\, : \,\,\,  f\in  N^1\cup\dots\cup N^\ell\big\}
\;\subseteq\; 2^{[n]\backslash (a\cup b)}.
\]

{\bf Notation.}
Let $\km{a}:=\{f\,:\; f\subseteq a\}$
be the full simplex on $a$, $\partial(\km{a}):=\{f\,:\; f\subsetneq a\}$,
and similarly define $\km{b},\partial(\km{b})$ from $b$.
If $Y$ and $Z$ are subcomplexes (or just subsets) of $X$ on disjoint
sets of vertices, then their {\em join} is
the simplicial complex
$\, Y\ast Z  \, := \, \{\,f \cup g
\, \,:\,\; f\in Y,\, g\in Z\}$.
In particular, $\,\{a\}\ast F = N^1\cup\dots\cup N^\ell$, and
it is straightforward to check that the triple join
$\,\km{a}\ast F \ast \partial(\km{b})\,$ is
a subcomplex of $X$.


\begin{theorem}
\klabel{flip}
The monomial ideal $M_X^\prime$
is square-free if and only if
$\,(\km{a}\ast F \ast \km{b}) \cap X =
\km{a}\ast F \ast \partial(\km{b})$
if and only if $\,\km{a}\ast F \ast \{b\}$ is disjoint from $X\,$
if and only if $\,X\cap (F \ast \{b\}) =\emptyset$.
If this holds then the neighboring simplicial complex $X^\prime$
with $M_{X'} = M_X^{\prime}$ is given
\vspace{-0.5ex}
by
\begin{equation}
\label{flipneighbor}
X^\prime
\quad = \quad
\big( X
\backslash (\km{a}\ast F \ast \partial(\km{b}))
\big)
\;\cup\; (\partial(\km{a})\ast F \ast \km{b})\,.
\vspace{1ex}
\end{equation}
\end{theorem}

Theorem \ref{flip} describes all the edges $\{X,X'\}$ in the {\em graph of
simplicial complexes}, that is, the subgraph of $\CG$ induced  on
square-free monomial ideals. The transition from $X$ to $X'$
generalizes the familiar notion of {\em bistellar flips}.
They correspond to the case $\ell=m=1$
and $b\notin X$. Here the condition in the first sentence
of Theorem \ref{flip} is automatically satisfied, meaning that
the  $c$-neighbor $M^\prime_X$ of $M_X$ is square-free.
These bistellar flips are a standard tool for locally altering
combinatorial manifolds (see e.g.~\cite{Viro}) or
triangulations of point configurations
(see e.g.~\cite{MaTh}).
\par


\begin{ex}
\klabel{bnotinX}
\rm
Let $a$ be an edge in a triangulated manifold $X$ of dimension two.
If $a\cup b$ supports the two triangles meeting along $a$,
then $N_c=\{a\}$, $\tilde{N}_c=\emptyset$, i.e.,
$\ell=m=1$, $N^1=N_c$, and $F=\{\emptyset\}$.
We are in the $b\notin X$ case,
and $\,\km{a}\ast F \ast \partial(\km{b})\,$ consists of the two
triangles and their faces.
\end{ex}

\hspace*{\fill}
\unitlength=1.0pt
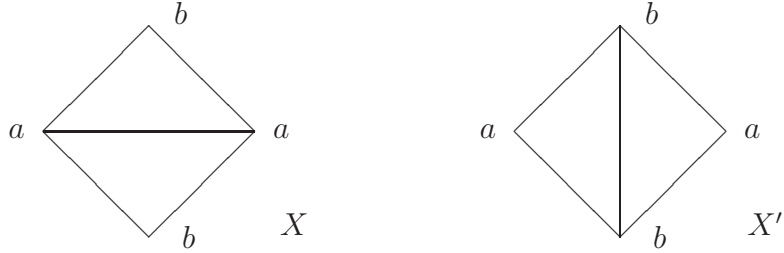
\begin{picture}(100.00,90.00)(10.00,40.00)
\put(52.00,125.00){\makebox(0.00,0.00){$b$}}
\put(55.00,40.00){\makebox(0.00,0.00){$b$}}
\put(90.00,80.00){\makebox(0.00,0.00){$a$}}
\put(-10.00,80.00){\makebox(0.00,0.00){$a$}}
\put(95.00,45.00){\makebox(0.00,0.00){$X$}}
\put(0.00,80.00){\line(1,0){80.00}}
\put(40.00,40.00){\line(-1,1){40.00}}
\put(80.00,80.00){\line(-1,-1){40.00}}
\put(40.00,120.00){\line(1,-1){40.00}}
\put(0.00,80.00){\line(1,1){40.00}}
\end{picture}
\hspace*{\fill}
\unitlength=1.0pt
\begin{picture}(100,90.00)(10.00,40.00)
\put(52.00,125.00){\makebox(0.00,0.00){$b$}}
\put(55.00,40.00){\makebox(0.00,0.00){$b$}}
\put(90.00,80.00){\makebox(0.00,0.00){$a$}}
\put(-10.00,80.00){\makebox(0.00,0.00){$a$}}
\put(95.00,45.00){\makebox(0.00,0.00){$X^\prime$}}
\put(40.00,40.00){\line(0,1){80.00}}
\put(40.00,40.00){\line(-1,1){40.00}}
\put(80.00,80.00){\line(-1,-1){40.00}}
\put(40.00,120.00){\line(1,-1){40.00}}
\put(0.00,80.00){\line(1,1){40.00}}
\end{picture}
\vspace{1ex}
\hspace*{\fill}
\par

>From $X$, we remove $\,\{a\}\ast F \ast \partial(\km{b})$,
i.e., the two triangles and their common edge.
They are replaced, in $X^\prime$, by the two triangles
$\,\partial(\km{a})\ast F \ast \{b\}\,$ with common edge $b$.
\par


\begin{ex}
\klabel{abinX}
\rm We still consider a triangulated surface.
Let $a$ be a trivalent vertex being adjacent to an edge $b$
and a third vertex $A$. In particular, $a\cup b\in X$.
\end{ex}

\hspace*{\fill}
\unitlength=1.0pt
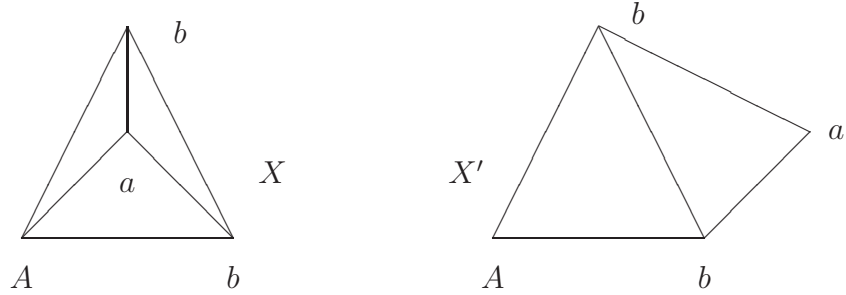
\begin{picture}(100,99.00)(0.00,0.00)
\put(60.00,98.00){\makebox(0.00,0.00){$b$}}
\put(40.00,40.00){\makebox(0.00,0.00){$a$}}
\put(80.00,5.00){\makebox(0.00,0.00){$b$}}
\put(0.00,5.00){\makebox(0.00,0.00){$A$}}
\put(95.00,45.00){\makebox(0.00,0.00){$X$}}
\put(40.00,60.00){\line(0,1){40.00}}
\put(40.00,100.00){\line(1,-2){40.00}}
\put(0.00,20.00){\line(1,2){40.00}}
\put(80.00,20.00){\line(-1,0){80.00}}
\put(40.00,60.00){\line(1,-1){40.00}}
\put(0.00,20.00){\line(1,1){40.00}}
\end{picture}
\hspace*{\fill}
\unitlength=1.0pt
\begin{picture}(100,99.00)(0.00,0.00)
\put(55.00,105.00){\makebox(0.00,0.00){$b$}}
\put(130.00,60.00){\makebox(0.00,0.00){$a$}}
\put(80.00,5.00){\makebox(0.00,0.00){$b$}}
\put(0.00,5.00){\makebox(0.00,0.00){$A$}}
\put(-10,45.00){\makebox(0.00,0.00){$X^\prime$}}
\put(40.00,100.00){\line(1,-2){40.00}}
\put(0.00,20.00){\line(1,2){40.00}}
\put(80.00,20.00){\line(-1,0){80.00}}
\put(80.00,20.00){\line(1,1){40.00}}
\put(40.00,100.00){\line(2,-1){80.00}}
\end{picture}
\hspace*{\fill}
\par

Here $N_c=\{\ko{Aa}\}$, $\,\tilde{N}_c=\emptyset$,
$\,\ell=m=1$, $\,N^1=N_c$, $\,F=\{A\}$. One obtains $X^\prime$
from $X$ by removing the edge $\ko{Aa}$ together with the adjacent triangles
and, afterwards, adding the triangle
formed by edge $b$ and vertex $A$. The new complex $X'$
is no longer part of a triangulation of a
two-dimensional manifold, since the edge $b$ is incident
to three triangles in $X'$.
The geometric realization $|X^\prime|$ looks like $|X|$
plus the additional triangle $\Delta(bba)$ sticking out of it.
\par



\kc{
\begin{ex}
\klabel{mis2}
\rm
Finally, we would like to show that
flipping backwards Example \ref{abinX} gives an instance with $m=2$,
i.e., with more than one neighbor in a fixed tangent direction $c$.
Let $X$ consist of three triangles
$\Delta(AST)$, $\Delta(BST)$, $\Delta(CST)$ sharing the common edge
$a:=\ko{ST}$.
With $b:=\{B\}$, we obtain
$\,N_c=\{\Delta(AST), \Delta(CST)\}$ and $\,\tilde{N}_c=\emptyset$.
Since $\ko{ST}$ does not belong to $N_c$, this yields
$\,m=2$, $\,N^1=\{\Delta(AST)\}$, and  $\,N^2=\{\Delta(CST)\}$.
Now, choosing $F$ among $\{A\}$, $\{B\}$, or $\{A,B\}$,
we have three possibilities to construct neighbors of $X$.
In each case, the square free condition of Theorem \ref{flip}
is satisfied, and we
obtain the following results for $X^\prime$:
\kitem
\item[(A)]
$X^\prime = \big\{\Delta(BTA), \Delta(BST), \Delta(BAS), \Delta(CST)\big\}
\cup\big\{\mbox{faces}\big\}$; this is like the $X$
from Example \ref{abinX}.
\item[(B)]
$X^\prime = \big\{\Delta(BTC), \Delta(BST), \Delta(BCS), \Delta(AST)\big\}
\cup\big\{\mbox{faces}\big\}$; this looks similar as the previous complex --
but now the other triangle $\Delta(CST)$ has been subdivided.
\item[(AB)]
$X^\prime = \big\{\Delta(BTA), \Delta(BST), \Delta(BAS), \Delta(BTC),
\Delta(BCS)\big\}\cup\big\{\mbox{faces}\big\}$,
i.e., both triangles $\Delta(AST)$ and $\Delta(CST)$
are subdivided by the same inner vertex $B$.
\vspace{2ex}
\kenditem
\end{ex}
\par
}

{\sl Proof of Theorem \ref{flip}: }
Recall from Corollary \ref{neighbors} that
$M_X^\prime$ has two types of generators.
The first were $\kx^u$ with $u\notin X$ and
$(u\cup a)\backslash b \notin N^{1}\cup\dots\cup N^{\ell}$,
the second were $\kx^{v+c}$ with $v\notin X$ and
$(v\cup a)\backslash b \in N^{1}\cup\dots\cup N^{\ell}$.
Hence, $M_X^\prime$ is square-free if
and only if every non reduced generator
$v+c$ of the latter type
finds some reduced generator $g$ with $g\leq v+c$.
Since this implies
$g\subseteq (v\cup a)\backslash b \in N^{1}\cup\dots\cup N^{\ell}
\subseteq N_c \subseteq X$, the generator $g$
cannot be of type one. If there is a type two generator
 $g=u+c\leq v+c$,
then, since $u\cap a =\emptyset$, we obtain $u\subseteq v\backslash a$,
and $u\notin X$ implies $\,v\backslash a \notin X$.
On the other hand, if
\kc{$\,v\backslash a \notin X$,}
then
$u:=v\backslash a$ indeed does the job.
We conclude that $M_X^\prime$ is square-free if and only if
there is no $v\in \Delta_{n-1}\backslash X$ with
$(v\cup a) \backslash b \in N^{1}\cup\dots\cup N^{\ell}$
and $\,v\backslash a \in X$. This condition is equivalent to the one
stated in Theorem \ref{flip}. To see this
take  $\,f:=v\backslash(a\cup b)\in F \,$ or $\,v :=f\cup(a\cup b)$.

Now, let us assume that this square-free condition is satisfied
and take equation (\ref{flipneighbor}) of the theorem as a definition
of some subset $X'\subset\Delta_{n-1}$.
We first show that $X'$ is indeed a
simplicial complex. Afterwards, we will
check that $M_{X^\prime}=M_X^\prime$.
\vspace{0.7ex}\\
{\sl Step 1}\,.
We claim the following:
{\sl
Let $f^\prime$ be a subset of $f$ which lies in $2^{[n]\backslash (a\cup b)}$.
Then, $\,f^\prime\in F\,$ and \,$ (a\cup f)\in X \,$
if and only if $\, f\in F\, $ and $\,(a\cup f^\prime\cup b)\notin X $.}\\
If $F$ was refering to $N_c$ instead to its subset
$N^1\cup\dots\cup N^\ell$, the claim would follow directly from the definition.
On the other hand, if both $\,a\cup f^\prime\,$ and
$\,a\cup f\,$ are in  $\, N_c$, then the whole flag
in between belongs to $N_c$ and, moreover, to the same connected component of
$\langle N_c \rangle$. In particular, $f^\prime \in F\,$ if and only $f\in F$.
\vspace{0.7ex}\\
{\sl Step 2}\,.
{\sl
$\,X^\prime \,= \, X
\backslash \big[\{a\}\ast F \ast \partial(\km{b})\big]
\;\cup\; \big[\partial(\km{a})\ast F \ast \{b\}\big]\,$
is a simplicial complex}:\\
First, we check that $\{a\}\ast F \ast \partial(\km{b})$
is, inside $X$, closed under enlargement.
Let $(a\cup f^\prime\cup b^{\prime\prime}) \subseteq
(a\cup f\cup b^{\prime})\in X$
with  $f^\prime\in F$, $b^{\prime\prime}\subsetneq b$,
$f\in 2^{[n]\backslash (a\cup b)}$, and
$b^{\prime} \supseteq b^{\prime\prime}$.
Then, since $\,(a\cup f)\in X$, we may use Step 1 to obtain $f\in F$.
Moreover, since
$\,X\cap [\km{a}\ast F \ast \{b\}]= \emptyset$,
the set $\,b^{\prime\prime}$ cannot equal $b$.\\
Now, we check the subsets of the  elements of $X^\prime\backslash X$.
Take, w.l.o.g.,
$g:=(a^{\prime\prime}\cup f^\prime\cup b)\subseteq (a^\prime \cup f\cup b)$
with $a^{\prime\prime}\subseteq a^\prime \subsetneq a$ and $f\in F$.
If $(a\cup f^\prime\cup b)\in X$, then $g\in X$, and we are done.
If not, then Step 1 implies that $f^\prime \in F$,
\vspace{0.7ex}
hence $g\in \big[\partial(\km{a})\ast F \ast \{b\}\big]\,$.\\
{\sl Step 3}\,.
{\sl $M_X^\prime= M_{X^\prime}$}:
Translating the square-free $M^\prime_X$ generators into the
$(\km{a}*F*\km{b})$-language, we obtain as exponents
$(a'\cup f\cup b')\notin X$ with $f\notin F$
for those of the first type and
$({a\cup f})$ with $f\in F$ for those of the second.
(The condition $(f\cup b)\notin X$ follows automatically from
$X\cap [\km{a}\ast F \ast \{b\}]=\emptyset$.)
While the type two generators are the minimal elements
of $[\{a\}\ast F \ast \partial(\km{b})\big]\subseteq X\backslash X'$,
the type one generators do neither belong to $X$, nor to
the part being changed during the transition to $X'$.
Hence, it remains to consider $g=(a^\prime\cup f\cup b^\prime)$
belonging to neither $X$, nor
$\big[\partial(\km{a})\ast F \ast \{b\}\big]$ and to show
that $\kx^g\in M_X^\prime$.
If $f\notin F$, then $g$ is obviously a generator of type one.
Assuming $f\in F$, then there are two possibilities:
If $a^\prime = a$, then $g\supseteq (a\cup f)$, i.e., a type two generator
takes care.
If $a^\prime \subsetneq a$, then $b^\prime \subsetneq b$,
hence $(a\cup f)\cup b^\prime\notin X$ implies
$(a\cup f) \in \tilde{N}_c$, and
we obtain a contradiction to $f\in F$.
\qed

\section{The Next Steps}\klabel{problems}

The following  list of open problems
arises naturally from our investigations.

{\bf Problem.}
Does the graph $ \CG_{n,\kk} $ depend on the field $\kk$?
For instance, is ${\cal G}_{n, {\R}}  =
{\cal G}_{n, {\C}} $~?
For $\kk$ algebraically closed, does
$\CG_{n,\kk}$ depend on the characteristic of $\kk$?
\par

While Theorem \ref{flip} shows that the graph of simplicial complexes
is independent of the field $\kk$,
the following example suggests that the answer might be ``yes'',
anyway.
Let $n = 2$, $M_1 = \langle x^4, y^2 \rangle$,
and $M_2 = \langle x^2, y^4 \rangle$.
The edge ideals connecting $M_1$ and $M_2$ have the form
\[
I \quad = \quad \langle
\, x^4, \,  y^2+a_1 y x+a_2 x^2 \,\rangle \quad = \quad
\langle \,x^2+b_1 y x+b_2 y^2, \, y^4 \, \rangle,
\]
where $a_1,a_2,b_1,b_2$ are scalars in $\kk$ satisfying
$\,
a_2 b_2 - 1  =
a_1 - a_2 b_1  =
 b_1^3 - 2 b_1 b_2  =  0
$.
These three equations define a scheme
which is the reduced
union of two irreducible components if $\mbox{\rm char}\,\kk\neq 2$,
and which is non-reduced but irreducible if $\mbox{\rm char}\,\kk= 2$.
\par

{\bf Problem.}
Find a purely combinatorial description for the graph of partitions.
What is the exact relationship with the directed graphs
studied by Evain in \cite{Eva}~?

{\bf Problem.}
Does the induced subgraph property hold for all
toric Hilbert schemes~?

{\bf Problem.}
Do there exist monomial ideals $M_1,M_2\subseteq\kk[\kx]$
having the same Hilbert function with respect to
two different gradings $\Z^n/\Z c$ and $\Z^n/\Z c'$?
Or is this
\kc{rather common}
provided that
$\Omega(M_1,M_2)=\emptyset$?
\par

Our results in Sections \ref{algorithm} and \ref{Hilbert}
provide a method for
constructing the graph of a multigraded Hilbert scheme $\Hilb_h$,
provided the following problem has been solved.

{\bf Problem.}
Develop a  practical algorithm for
computing all monomial ideals in $\Hilb_h$.
\par

Peeva and Stillman \cite{PeSt} recently proved a connectivity
theorem for Hilbert schemes over the exterior algebra,
using methods similar to those in Section \ref{simplicial}.
Monomial ideals in the exterior algebra
being square-free, the following question arises.

{\bf Problem.}
Given any grading and Hilbert function on the
exterior algebra, is its flip graph the same
as the subgraph induced from
our graph of simplicial complexes?

It is natural to wonder about the topological significance
of these flips.

{\bf Problem.}
Which topological invariants remain unchained by the
generalized flips of Theorem \ref{flip}. How can one decide,
by means of a practical algorithm, whether two given simplicial complexes
can be connected by a chain of those flips?
\par

\medskip

{\sl Acknowledgement.} We would like to thank Rekha Thomas
for carefully reading a manuscript version of this paper.

%
%

\vfill

{\small
\setbox0\hbox{FB Mathematik und Informatik, WE2}
\parbox{\wd0}{
Klaus Altmann\\
FB Mathematik und Informatik, WE2\\
Freie Universit\"at Berlin\\
Arnimallee 3\\
D-14195 Berlin, Germany\\
email: altmann@math.fu-berlin.de}
\setbox1\hbox{University of California at Berkeley}\hfill\parbox{\wd1}{
Bernd Sturmfels\\
Department of Mathematics\\
University of California at Berkeley\\
Berkeley, CA 94720, USA\\
email: bernd@math.berkeley.edu}}


{\small
\begin{thebibliography}{A}

\bibitem[1]{AlCh} K.~Altmann and J.A.~Christophersen:
Deforming Stanley-Reisner rings,   \hfill \break
{\tt  math.AG/0006139}.

\bibitem[2]{Eis} D.~Eisenbud: {\em Commutative algebra. With a
view toward algebraic geometry}. Graduate Texts in Mathematics,
150. Springer-Verlag, New York, 1995.

\bibitem[3]{Eva} L.~Evain:
Incidence relations among the Schubert cells of equivariant Hilbert schemes,
to appear in {\sl Mathematische Zeitschrift}, {\tt math.AG/0005233}.


\bibitem[4]{Ha} R.~Hartshorne: Algebraic Geometry. Graduate Texts in
Mathematics {\bf 52}, Springer-Verlag 1977.

\bibitem[5]{HaSt} M.~Haiman and B.~Sturmfels:  Multigraded Hilbert schemes,
{\tt math.AG/0201271}.

\bibitem[6]{MaTh} D.~Maclagan and R.~Thomas:
Combinatorics of the toric Hilbert scheme,
{\em Discrete Comput.~Geom.} {\bf 27} (2002) 249--272.

\bibitem[7]{PeSt} I.~Peeva and M.~Stillman:
Exterior Hilbert schemes, in preparation.

\bibitem[8]{Ree}A.~Reeves: The radius of the Hilbert scheme,
{\em  J. Algebraic Geom.} {\bf 4} (1995) 639--657.

\bibitem[9]{San} F.~Santos:
Non-connected toric Hilbert schemes,
preprint, {\tt  math.CO/0204044}.

\bibitem[10]{Viro} O.~Viro:
Lectures on combinatorial presentations of manifolds,
Differential geometry and topology (Alghero, 1992),
World Sci.\ Publishing, 1993, 244--264.


\end{thebibliography}
}
\end{document}